\documentclass[11pt,a4paper]{article}
\usepackage[top=2.5cm,bottom=2.5cm,left=2.2cm,right=2.2cm]{geometry}
\usepackage{tikz}
\usetikzlibrary{positioning}
\usepackage{lineno}
\usepackage[T1]{fontenc}
\usepackage[utf8]{inputenc}
\usepackage{color}
\usepackage{graphicx}
\usepackage{amsfonts}
\usepackage{extarrows}
\usepackage{amsmath,amsthm,amssymb,color}
\usepackage{hyperref}
\usepackage{eepic}
\usepackage{lineno}
\usepackage{enumerate}	
\usepackage{paralist}
\usepackage{cite}
\usepackage{algorithm}
\usepackage{algorithmicx}
\usepackage{algpseudocode}
\usepackage{authblk}
\usepackage{mathtools}
\usepackage{pifont}
\usepackage[numbers,sort&compress]{natbib}

\usepackage{tikz}

\hypersetup
{
	colorlinks=true,
	linkcolor=blue,
	filecolor=blue,
	urlcolor=blue,
	citecolor=cyan,
}

\newtheorem{theorem}{Theorem}[section]

\newtheorem{corollary}[theorem]{Corollary}
\newtheorem{conjecture}[theorem]{Conjecture}

\theoremstyle{definition}

\newtheorem{claim}{\indent Claim}
\newtheorem{remark}[theorem]{Remark}

\newtheorem{case}{\indent Case}[section]

{\setlength{\leftmargini}{1.5\parindent}
 \begin{itemize}
 \setlength{\itemsep}{-1.1mm}}
{\end{itemize}}

\begin{document}
	\title{\bf Every $2k$-connected $(P_2\cup kP_1)$-free graph with toughness greater than one  is hamiltonian-connected}
	
	\author{\bf Feng Liu\footnote{Email: liufeng0609@126.com.}}
	\affil{Department of Mathematics,
		East China Normal University, Shanghai, 200241, China}
	\date{}
	\maketitle
\begin{abstract}
 Given a graph $H$, a graph $G$ is $H$-free if $G$ does not contain $H$ as an induced subgraph.  Shi and Shan conjectured that every $1$-tough $2k$-connected $(P_2 \cup kP_1)$-free graph is hamiltonian for $k \geq 4$. This conjecture has been independently confirmed by Xu, Li, and Zhou, as well as by Ota and Sanka.
  Inspired by this, we prove that every $2k$-connected  $(P_2\cup kP_1)$-free graph with toughness greater than one is hamiltonian-connected.

		
\end{abstract}
	
\section{Introduction}
We consider finite simple  graphs, and use standard terminology and notations from \cite{Bondy2008, West1996} throughout this article. We denote by $V(G)$ and $E(G)$ the vertex set and edge set of a graph $G,$ respectively, and denote by $|G|$ and $e(G)$ the order and size of $G,$ respectively. For a vertex $x\in V(G)$ and a subgraph $H$ of $G$, $N_H(x)$ denotes  the set of neighbors of $x$ that are contained in $V(H)$. For a vertex subset $S\subseteq V(G)$, define $N_G(S)=\cup_{x\in S}N_G(x)\setminus S$ and $N_H(S)=N_G(S)\cap V(H)$. If $F$ is a subgraph of $G$, we write $N_F(H)$ for $N_F(V(H))$. We use $G[S]$ to denote the subgraph of $G$ induced by $S$, and let $G-S=G[V(G)\setminus S]$.  Furthermore, we denote the number of components in a graph $G$ by  $c(G)$. The symbol $[k]$ used in this article represents the set $\{1,2,\ldots,k\}$. The subscript $G$ will be omitted in all the notation above if no confusion may arise. 

For two distinct vertices $x$ and $y$ in $G$, an $(x, y)$-path is a path whose endpoints are $x$ and $y$.  Let $P$ be a path. We   use $P[u,v]$ to denote the subpath of $P$ between two  vertices $u$  and $v$. 
A Hamilton path in $G$ is a path containing every vertex of $G$.  A Hamilton  cycle in $G$ is a cycle containing every vertex of $G$. The graph $G$ is hamiltonian if it contains a Hamilton cycle. 
A graph is called  hamiltonian-connected if between any 
two distinct vertices there is a Hamilton path. 

For a given graph $H$, a graph $G$ is called $H$-free if $G$ does not contain $H$ as an induced subgraph. For vertex disjoint graphs $H$ and $F$, $H\cup F$ denotes the disjoint union of graphs $H$ and $F$. A linear forest is a graph consisting of disjoint paths.
As usual, $P_n$ denotes the path on $n$ vertices. For positive integer $k$ and $\ell$, $kP_{\ell}$ denotes the linear forest consisting of $k$ disjoint copies of the path $P_{\ell}$.

For a positive integer $k$, a connected graph $G$ is said to be $k$-connected if any deletion of at most $k-1$ vertices on $G$ also results in a connected graph. For a graph $G$ with $S\subset V(G)$, denote by $G[S]$ the subgraph of $G$ induced by $S$. The toughness of a graph $G$, denoted by $\tau(G)$, is defined as
\begin{flalign*}
    \tau(G)=\min \{\frac{|S|}{c(G-S)}: S\subseteq V(G), c(G-S)\ge 2\}
\end{flalign*}
if $G$ is not a complete graph and $\tau(G)=\infty$ otherwise.
For a positive real number $t$,
a graph $G$ is called $t$-tough if $\tau(G)\ge t$, that is, $|S|\ge t\cdot c(G-S)$ for each $S\subseteq V(G)$ with $c(G-S)\ge 2$.

The concept of toughness of a graph was introduced by Chv\'{a}tal \cite{Ch}.
Clearly, every  hamiltonian graph is $1$-tough, but the converse is not true.
Chv\'{a}tal \cite{Ch} proposed the following conjecture, which is known as Chv\'{a}tal's toughness conjecture.

\begin{conjecture}[Chv\'{a}tal \cite{Ch}]\label{Or}
There exists a constant $t_0$ such that every $t_0$-tough graph
with at least three vertices is hamiltonian.
\end{conjecture}

Bauer, Broersma and Veldman \cite{BBV} showed that $t_0\ge \frac{9}{4}$ if it exists.
Conjecture \ref{Or}  has been confirmed  for a number of special classes of graphs
\cite{BBS,BHT,BPP,BXY,CJKJ,DKS,GS,HG,KK,LBZ,OS,Shan0,Sh,Shan,SS,Wei}.
For example, it has been confirmed for graphs with forbidden (small) linear forests, such as $1$-tough $R$-free graphs with $R=P_3\cup P_1, P_2\cup 2P_1$ \cite{LBZ},
$2$-tough $2P_2$-free graphs \cite{BPP,Sh,OS}, $3$-tough $(P_2\cup 3P_1)$-free graphs  \cite{HG}, $7$-tough $(P_3\cup 2P_1)$-free graphs \cite{GS} and $15$-tough $(P_3\cup P_2)$-free graphs \cite{Shan}. Though great efforts have been made,  it remains open. In 2022, Shi and Shan \cite{SS}  proposed the following conjecture.

\begin{conjecture}[Shi-Shan \cite{SS}]\label{Shi}
	 Let $k\ge 4$ be an integer and let $G$ be a
	$1$-tough and $2k$-connected $(P_2\cup kP_1)$-free graph. Then $G$ is hamiltonian.
\end{conjecture} 
Recently, Conjecture~\ref{Shi} has been independently confirmed by Xu, Li, and Zhou, as well as by Ota and Sanka.
It is well known that a hamiltonian-connected graph has toughness strictly larger than one. 
Suppose $R$ is a graph such that
every $1$-tough $R$-free graph is hamiltonian. Is then every $R$-free graph $G$ with $\tau(G)>1$ hamiltonian-connected?
In 1978, Jung \cite{Jung} obtained the following result, in which he showed that for
$P_4$-free graphs, the necessary condition $\tau(G)>1$ is also a sufficient condition for hamiltonian-connectivity.

\begin{theorem}[Jung \cite{Jung}]
	Let $G$ be a $P_4$-free graph. Then $G$ is hamiltonian-connected if and only if $\tau(G)>1$.
\end{theorem}

 In 2022, Zheng, Broersma and Wang \cite{Zheng} provided some sufficient conditions for hamiltonian-connectivity by replacing other graphs to $P_4$.

\begin{theorem}[Zheng-Broersma-Wang \cite{Zheng}]
	Let $R$ be an induced subgraph of $K_1\cup P_3$ or $2K_1\cup K_2$. Then every
	$R$-free graph $G$ with $\tau(G)>1$ on at least three vertices is hamiltonian-connected.
\end{theorem}

There is a classic result about hamiltonian-connected, due to Chv\'{a}tal and Erd\H{o}s \cite{CE}.

\begin{theorem}[Chv\'{a}tal-Erd\H{o}s \cite{CE}]\label{coe}
	For any integer $k\ge 1$, every $k$-connected $kP_1$-free graphs on at least three vertices is hamiltonian-connected.
\end{theorem}

Note that a $k$-connected $kP_1$-free graph must have toughness greater than one. However, the constant connectivity condition cannot guarantee the property of hamiltonian-connected in $(P_2\cup kP_1)$-free graphs.
In this paper, we show the following result.

\begin{theorem}\label{MTHM-1}
Let $k$ be a positive integer. Every $2k$-connected  $(P_2\cup kP_1)$-free graph with toughness greater than one is hamiltonian-connected.
\end{theorem}
Since  a non-complete $k$-tough graph is $\lceil 2k \rceil$-connected, as a corollary of  Theorem \ref{MTHM-1}, we get the following result.

\begin{corollary} \label{zhong}
For any positive integer $k\geq 2$, every $k$-tough $(P_2\cup kP_1)$-free graph is hamiltonian-connected.
\end{corollary}
\section{Proof of Theorem \ref{MTHM-1}}
The aim of this section is to prove Theorem \ref{MTHM-1}. 
Before our proof we list some useful notations needed in later proofs. 
Let $P$ be an oriented path. For $u\in V(P)$,
denote by $u^{+1}$ the immediate successor of $u$ and $u^{-1}$ the immediate predecessor of $u$ on $P$.
For an integer  $\ell\ge 2$,
denote by $u^{+\ell}$ the immediate successor of $u^{+(\ell-1)}$ and $u^{-\ell}$ the immediate predecessor of $u^{-(\ell-1)}$ on $P$.
For convenience, we write  $u^{+}$ for  $u^{+1}$ and $u^{-}$ for  $u^{-1}$.
Let $x$ be the endpoint of $P$ such that $x^+$ does not exist.
For $S\subseteq V(P)$, let $S^+=\{u^+:u\in S\setminus \{x\}\}$. 
For $u,v\in V(P)$,  $\overrightarrow{P}[u,v]$ denotes the segment of $P$ from $u$ to $v$ which follows the orientation of $P$, while $\overleftarrow{P}[u,v]$ denotes the opposite segment of $P$ from $v$ to $u$. Particularly, if $u=v$, then $\overrightarrow{P}[u,v]=\overleftarrow{P}[u,v]=u$.

We are now ready to prove our main result.
\begin{proof}[\bf Proof of Theorem \ref{MTHM-1}.]

 Let $u,v\in V(G)$ be any two distinct vertices and $P$ be a longest $(u,v)$-path in $G$.  We may assume that $P$ is not a Hamilton path of $G$. Let $H$ be a component of $G-V(P)$ and let $N_P(H)=\{x_1,x_2,\ldots, x_t\}$. 
 \begin{claim}\label{Claim-independent}
    $N_P(H)^+$ is an independent set.
 \end{claim}
  Otherwise, suppose  $x_i^+$ and $x_j^+$ are adjacent for some $i,j\in [t]$. Furthermore, assume by symmetry that $i<j$. Since $x_i,x_j$ both have a neighbor in $H$, there is a path $Q$ between $x_i,x_j$ with the interior in $H$. Consequently,
  \begin{flalign*}
      u,\overrightarrow{P}[u,x_i],x_i,Q,x_j,\overleftarrow{P}[x_j,x_i^+],x_i^+,x_j^+,\overrightarrow{P}[x_j^+,v],v
  \end{flalign*}
is a $(u,v)$-path longer than $P$, a contradiction.	This proves Claim \ref{Claim-independent}.

Since $P$ is a longest $(u,v)$-path, $H$ has no two consecutive neighbors on $P$. Thus, $x_ix_{i+1}\notin E(P)$ for any $i\in [t-1]$. That is, $N_P(H)^+\cap N_P(H)=\emptyset$.

\begin{claim}\label{Claim-isolated}
$H$ is an isolated vertex.
\end{claim}
Since $G$ is $2k$-connected, $|N_P(H)^+|\geq |N_P(H)|-1\geq 2k-1$. Suppose to the contrary that $H$ contains an edge $zw$. By Claim \ref{Claim-independent}, we have that $G[V(N_P(H)^+)\cup \{z,w\}]$ contains  $P_2\cup (2k-1)P_1$ as a induced subgraph. This implies that $G$ contains $P_2\cup kP_1$ as an induced subgraph, a contradiction. This proves Claim \ref{Claim-isolated}.

By Claim \ref{Claim-isolated}, we have that $H$ is an isolated vertex. For convenience, we may assume that $V(H)=\{x\}$. For $i\in [t-1]$, denote $S_i=V(P[x_i^+,x_{i+1}^-])$.  Clearly, $|S_i|\geq 1$ for $i\in[t-1]$. In addition, denote $S_0=V(P[u,x_{1}^-])$ and $S_t=V(P[x_t^+,v])$. Furthermore, for $i\in [t-1]$, let $X\subseteq N_P(x)^+$  contain $x_i^+$ with $|X|=2k-1$.

\begin{claim}\label{Claim-|X|=2k-1}
For $i\in [t-1]$, if $j$ is an integer with $1\leq j\leq |S_i|$, then $|N_P(x_i^{+j})\cap X|=0$ if $j$ is odd integer, and
$|N_P(x_i^{+j})\cap X|\geq k+1$ if $j$ is even integer.
\end{claim}

We use induction on $j$. It is not difficult to verify that Claim \ref{Claim-|X|=2k-1} holds for $j=1$. Suppose $j=2$. If $|N_P(x_i^{+2})\cap X|\leq k$, then $|(X\setminus N_P(x_i^{+2}))\cup \{x,x_i^+, x_i^{+2}\}|\geq k+2$, and the edge $x_i^+x_i^{+2}$ is  the only edge in $G[(X\setminus N_P(x_i^{+2}))\cup \{x,x_i^+, x_i^{+2}\}]$. This implies that $G$ contains $P_2\cup kP_1$ as an induced subgraph, a contradiction.
Next let $j\geq 3$ and we suppose that Claim \ref{Claim-|X|=2k-1} holds for  any  $p$ with $1\leq p\leq j-1$. We distinguish two cases according as $j$ is even or odd. 

\begin{case}
$j$ is odd.
\end{case}
By the inductive hypothesis, $|N_P(x_i^{+(j-1)})\cap X|\geq k+1$. Suppose $N_P(x^{+j})\cap X\neq \emptyset$. Let $x_r^+\in N_P(x^{+j})\cap X$. If $|N_P(x^{+j})\cap X|\leq k$, then $|(X\setminus N_P(x^{+j})\cup \{x, x_i^{+j},x_r^{+}\}|\geq k+2$, and the edge $x_r^+x_i^{+j}$ is  the only edge in $G[(X\setminus N_P(x^{+j})\cup \{x, x_i^{+j},x_r^{+}\}]$. This implies that $G$ contains $P_2\cup kP_1$ as an induced subgraph, a contradiction. Therefore, $|N_P(x^{+j})\cap X|\geq k+1$. However, it follows that 
\begin{flalign*}
|N_P(x_i^{+(j-1)})\cap N_P(x_i^{+j}) \cap X|
&\ge |N_P(x_i^{+(j-1)})\cap X|+|N_P(x_i^{+j}) \cap X|-|X|\\
&\ge  k+1+k+1-(2k-1)>2.
\end{flalign*}
Let $x_p^+,x_q^+\in N_P(x_i^{+(j-1)})\cap N_P(x_i^{+j})\cap X$ with $1\le p<q\le t$.

If $p\geq i+1$, then 
\begin{flalign*}
u,\overrightarrow{P}[u,x_i^{+(j-1)}],x_i^{+(j-1)},x_p^+,\overrightarrow{P}[x_p^+,x_q],x_q,x,x_p\overleftarrow{P}[x_p,x_i^{+j}],x_i^{+j},x_q^{+},\overrightarrow{P}[x_q^{+},v],v
\end{flalign*}
 is a $(u,v)$-path longer than $P$, a contradiction.

If $q\leq i$, then 

\begin{flalign*}
    u,\overrightarrow{P}[u,x_p],x_p,x,x_q,\overleftarrow{P}[x_q,x_p^+], x_p^+, x_i^{+(j-1)},\overleftarrow{P}[x_i^{+(j-1)},x_q^{+}], x_q^+,x_i^{+j},\overrightarrow{P}[x_i^{+j},v],v
\end{flalign*}
 is a $(u,v)$-path longer than $P$, a contradiction.

If $p\leq i< q$, then 

\begin{flalign*}
u,\overrightarrow{P}[u,x_p], x_p,x,x_q,\overleftarrow{P}[x_q,x_i^{+j}],x_i^{+j},x_p^+,\overrightarrow{P}[x_p^+,x_i^{+(j-1)}],x_i^{+(j-1)},x_q^{+},\overrightarrow{P}[x_q^{+},v],v
\end{flalign*}
is a $(u,v)$-path longer than $P$, a contradiction. This proves Claim \ref{Claim-|X|=2k-1} for $j$ is odd.

\begin{case}
$j$ is even.
\end{case}
By the inductive hypothesis,  $N_P(x_i^{+(j-1)})\cap X=\emptyset$.
If $|N_P(x_i^{+j})\cap X|\le k$, then $$| (X\setminus N_P(x_i^{+j}))\cup \{x,x_i^{+j},x_i^{+(j-1)}\}|\ge k+2,$$and the edge  $x_i^{+j}x_i^{+(j-1)}$ is the only edge in 
$G[(X\setminus N_P(x_i^{+(j-1)}))\cup \{x,x_i^{+j},x_i^{+(j-1)}\}]$. This implies that $G$ contains $P_2\cup kP_1$ as an induced subgraph, a contradiction.

Therefore, $\left|N_P(x_i^{+j})\cap X\right|\ge k+1$. This  proves Claim \ref{Claim-|X|=2k-1}.

\begin{claim}\label{Claim-emptyset}
For $i\in [t-1]$, if $j$ is an odd integer with $1\leq j\leq |S_i|$, then $N_P(x_i^{+j})\cap N_P(x)^+=\emptyset.$
\end{claim}
Let $\mathcal{F}=\{X\subseteq N_P(x)^+:~X~ \text{contains}~x_i^+~  \text{with}~|X|=2k-1\}$.
 By Claim \ref{Claim-|X|=2k-1}, we have that $N_P(x_i^{+j})\cap X=\emptyset.$ Therefore,
 \begin{flalign*}
    N_P(x_i^{+j})\cap N_P(x)^+=\bigcup\limits_{X\in \mathcal{F}}( N_P(x_i^{+j})\cap X)=\emptyset.
 \end{flalign*}
This proves Claim \ref{Claim-emptyset}.

By a similar analysis as in the proofs of Claim \ref{Claim-|X|=2k-1} and Claim \ref{Claim-emptyset}, we obtain that if $|S_t| \geq 1$  and $j$ is an odd integer with $1 \leq j \leq |S_t| $, then  $N_P(x_t^{+j}) \cap N_P(x)^+ = \emptyset$. Furthermore, considering the opposite direction of path $P$ from $v$ to $u$, we obtain that if $|S_0| \geq 1 $ and $j$ is an odd integer with $1 \leq j \leq |S_0|$, then  $N_P(x_1^{-j}) \cap N_P(x)^+ = \emptyset$.

\begin{claim}\label{Claim-S_i}
For $i\in [t-1]$, $|S_i|$ is odd.
\end{claim}
Otherwise, there exists $i\in [t-1]$ such that $|S_i|$ is even. By Claim \ref{Claim-|X|=2k-1}, we have that $|N_P(x_i^{+|S_i|})\cap X|\geq k+1$. If $|N_P(x_{i+1})\cap X|\leq k-1$, then $|(X\setminus N_P(x_{i+1}))\cup \{x,x_{i+1}\} |\geq k+2$,  and the edge $xx_{i+1}$ is the only edge in $G[(X\setminus N_P(x_{i+1}))\cup \{x,x_{i+1}\}]$. This implies that $G$ contains $P_2\cup kP_1$ as an induced subgraph, a contradiction. Therefore, $|N_P(x_{i+1})\cap X|\geq k$.  However, it follows that 
\begin{flalign*}
|N_P(x_{i+1})\cap N_P(x_i^{+|S_i|}) \cap X|
&\ge |N_P(x_{i+1})\cap X|+|N_P(x_i^{+|S_i|}) \cap X|-|X|\\
&\ge  k+k+1-(2k-1)\geq 2.
\end{flalign*}

Let $x_p^+,x_q^+\in N_P(x_{i+1})\cap N_P(x_i^{+|S_i|}) \cap X$ with $1\le p<q\le t$.

If $p\geq i+1$, then 
\begin{flalign*}
u,\overrightarrow{P}[u,x_i^{+|S_i|}],x_i^{+|S_i|},x_p^+,\overrightarrow{P}[x_p^+,x_q],x_q,x,x_p\overleftarrow{P}[x_p,x_{i+1}],x_{i+1},x_q^{+},\overrightarrow{P}[x_q^{+},v],v
\end{flalign*}
 is a $(u,v)$-path longer than $P$, a contradiction.

If $q\leq i$, then 

\begin{flalign*}
    u,\overrightarrow{P}[u,x_p],x_p,x,x_q,\overleftarrow{P}[x_q,x_p^+], x_p^+, x_i^{+|S_i|},\overleftarrow{P}[x_i^{+|S_i|},x_q^{+}], x_q^+,x_{i+1},\overrightarrow{P}[x_{i+1},v],v
\end{flalign*}
 is a $(u,v)$-path longer than $P$, a contradiction.

If $p\leq i< q$, then 
\begin{flalign*}
u,\overrightarrow{P}[u,x_p], x_p,x,x_q,\overleftarrow{P}[x_q,x_{i+1}],x_{i+1},x_p^+,\overrightarrow{P}[x_p^+,x_i^{+|S_i|}],x_i^{+|S_i|},x_q^{+},\overrightarrow{P}[x_q^{+},v],v
\end{flalign*}
is a $(u,v)$-path longer than $P$, a contradiction. This proves Claim \ref{Claim-S_i}.

For $i\in [t]$, denote  $S_i^{\prime}=\{x_i^{+j}\in S_i:~ j ~\text{is odd}\}$. Furthermore, denote $S_0^{\prime}=\{x_1^{-j}\in S_0:~ j ~\text{is odd}\}$. Let $S'=\bigcup\limits_{i=0}^tS_i^{\prime}$.

\begin{claim}\label{Claim-S'}
$S^{\prime}$ is an independent set.
\end{claim}
Otherwise, we may assume that there is an edge $zw$ in $G[S']$.  Note that $N_P(z)\cap  N_P(x)^+=\emptyset$ and $N_P(w)\cap  N_P(x)^+=\emptyset$. By Claim \ref{Claim-independent}, we have $N_P(x)^+$ is an independent set. Therefore, the edge $zw$ is the only edge in $G[N_P(x)^+]\cup \{z,w\}$. Note that $|N_P(x)^+|\geq 2k-1$. This implies that $G$ contains $P_2\cup kP_1$ as an induced subgraph, a contradiction. This proves Claim \ref{Claim-S'}.

Note that $|S_i|$ is odd for $i\in [t-1]$. Clearly, $|V(P)|\leq 2|S^\prime|+1$. By Claim \ref{Claim-isolated}, we have that $V(G)\setminus V(P) $  is an independent set.

\begin{claim}\label{Claim-y}
For any $y\in V(G)\setminus V(P)$, $N_P(y)\cap S'=\emptyset$.
\end{claim}
 Clearly, $y\neq x$. Firstly, we prove that $ N_P(y)\cap N_P(x)^+=\emptyset$. Suppose $|N_P(y)\cap N_P(x)^+|\geq 2$. Without loss of generality, assume that $x_p^+,x_q^+\in N_P(y)$ with  $1\le p<q\le t$. However, it follows that  
\begin{flalign*}
    u,\overrightarrow{P}[u,x_p],x_p,x,x_q,\overleftarrow{P}[x_q,x_p^+],x_p^+,y,x_p^+,\overrightarrow{P}[x_p^+,v],v
\end{flalign*}
 is a $(u,v)$-path longer that $P$, a contradiction. Therefore, $|N_P(y)\cap N_P(x)^+|\leq 1$.  Note that $\{x\}\cup N_P(x)^+$ is an independent set. If $|N_P(y)\cap N_P(x)^+|=1$, then $G[N_P(x)^+\cup \{x,y\}]$ has order at least $2k+1$ and size one. This implies that $G[N_P(x)^+\cup \{x,y\}]$ contains $P_2\cup kP_1$ as an induced subgraph, and hence $G$ contains $P_2\cup kP_1$ as an induced subgraph, a contradiction. Therefore, $N_P(y)\cap N_P(x)^+=\emptyset$. 

 Now, we prove that $N_P(y)\cap S'=\emptyset$. Suppose to the contrary that there exists $z\in N_P(y)\cap S'$. Note that $z\notin N_P(x)^+$ and $N_P(z)\cap N_P(x)^+=\emptyset$.  Now, $G[N_P(x)^+\cup \{y,z\}]$ has order $2k+1$ and size one. This implies that $G[N_P(x)^+\cup \{y,z\}]$ contains $P_2\cup kP_1$ as an induced subgraph, and hence $G$ contains $P_2\cup kP_1$ as an induced subgraph, a contradiction. Therefore, $N_P(y)\cap S'=\emptyset$. This proves Claim \ref{Claim-y}.

Let $S^*=V(P)\setminus S'$. Then $|S^*|\leq \frac{|V(P)|+1}{2}$. Let $I=V(G)\setminus V(P)$. By Claim \ref{Claim-independent}, Claim \ref{Claim-S'} and Claim \ref{Claim-y}, we have that $I\cup S'$ is an independent set. Note that $|I|+|S'|\geq 1+\frac{|V(P)|-1}{2}\geq |S^*|$. However, it follows that 
\begin{flalign*}
\tau(G)\leq \frac{|S^*|}{c(G-S^*)}=\frac{|S^*|}{|I|+|S^\prime|}\leq \frac{|S^*|}{|S^*|} =1,
\end{flalign*}
contradicting to the condition that $\tau(G)>1$. This completes the proof of Theorem \ref{MTHM-1}. 
\end{proof}
\begin{remark}
    In this paper, we prove that every $2k$-connected  $(P_2\cup kP_1)$-free graph with toughness greater than one is hamiltonian-connected, where $k$ is a positive integer. Let $K_{s,t}$ denote the complete bipartite graph whose partite sets have cardinality $s$ and $t$, respectively.
    In our result, the condition of toughness greater than one is necessary. Consider the balanced complete bipartite graph $K_{\frac{n}{2},\frac{n}{2}}$ on $n$ vertices where $n\geq 4$ is even. It is easy to check $K_{\frac{n}{2},\frac{n}{2}}$ is a $1$-tough $\frac{n}{2}$-connected $(P_2\cup \frac{n}{4}P_1)$-free graph that is not hamiltonian-connected. 
\end{remark}

\section*{Acknowledgement} The author is grateful to Professor Xingzhi Zhan for his constant support and guidance. The author also thanks Chengli Li for careful reading of an early draft of this paper and for helpful comments. This research  was supported by the NSFC grant 12271170 and Science and Technology Commission of Shanghai Municipality (STCSM) grant 22DZ2229014.

\section*{Declaration}
\noindent$\textbf{Conflict~of~interest}$
The authors declare that they have no known competing financial interests or personal relationships that could have appeared to influence the work reported in this paper.
	
\noindent$\textbf{Data~availability}$
Data sharing not applicable to this paper as no datasets were generated or analysed during the current study.


\begin{thebibliography}{99}

\bibitem{BBS} D. Bauer, H.J. Broersma and E. Schmeichel,
Toughness in graphs--a survey,
\emph{Graphs Combin.}, 22 (2006) 1--35.

\bibitem{BBV} D. Bauer, H.J. Broersma and H.J. Veldman,
Not every $2$-tough graph is Hamiltonian,
\emph{Discrete Appl. Math.}, 99 (2000) 317--321.

\bibitem{BHT} T. B\"{o}hme, J. Harant and M. Tk\'{a}\v{c},
More than one tough chordal planar graphs are Hamiltonian,
\emph{J. Graph Theory}, 32  (1999) 405--410.

\bibitem{Bondy2008} J.A. Bondy and U.S.R. Murty, Graph theory. Graduate texts in mathematics, vol. 244. Springer, New York, 2008, pp. xii+651.	

\bibitem{BPP} H. Broersma, V. Patel and A. Pyatkin,
On toughness and hamiltonicity of $2K_2$-free graphs,
\emph{J. Graph Theory}, 75 (2014) 244--255.

\bibitem{BXY} H.J. Broersma, L. Xiong and K. Yoshimoto,
Toughness and hamiltonicity in $k$-trees,
\emph{Discrete Math.}, 307 (2007) 832--838.


\bibitem{CJKJ} G. Chen, M.S. Jacobson, A.E. K\'{e}zdy and L. Jen\"{o},
Tough enough chordal graphs are Hamiltonian,
\emph{Networks}, 31 (1998) 29--38.

\bibitem{Ch} V. Chv\'{a}tal, Tough graphs and Hamiltonian circuits,
\emph{Discrete Math.}, 5 (1973) 215--228.

\bibitem{CE} V. Chv\'{a}tal and P. Erd\"{o}s,
A note on Hamiltonian circuits, \emph{Discrete Math.}, 2 (1972) 111--113.

\bibitem{DKS} J.S. Deogun, D. Kratsch and G. Steiner,
$1$-tough cocomparability graphs are Hamiltonian,
\emph{Discrete Math.}, 170 (1997) 99--106.


\bibitem{GS} Y. Gao and S. Shan,
Hamiltonian cycles in $7$-tough $(P_3\cup 2P_1)$-free graphs,
\emph{Discrete Math.}, 345 (2022) 113069.

\bibitem{Jung} H.A. Jung, On a class of posets and the corresponding comparability graphs, \emph{J. Combin. Theory Ser. B}, 24 (1978) 125–133.

\bibitem{HG} A. Hatfield and E. Grimm,
Hamiltonicity of $3$-tough $(P_2\cup 3P_1)$-free graphs, arXiv: 2106.07083.

\bibitem{KK} A. Kabela and T. Kaiser, $10$-tough chordal graphs are Hamiltonian,
\emph{J. Combin. Theory Ser. B}, 122 (2017) 417--427.

\bibitem{LBZ} B. Li, H.J. Broersma and S. Zhang,
Forbidden subgraphs for hamiltonicity of $1$-tough graphs,
\emph{Discuss. Math. Graph Theory}, 36 (2016) 915--929.

\bibitem{OS} K. Ota and M. Sanka,
Hamiltonian cycles in $2$-tough $2K_2$-free graphs,
\emph{J. Graph Theory},
101 (2022) 769--781.

\bibitem{OS2} K. Ota and M. Sanka. Some conditions for hamiltonian cycles in $1$-tough $(K_2\cup kK_1)$-free graphs, \emph{Discrete Math.}, 347(3) (2024) 113841.


\bibitem{Shan0} S. Shan, An ore-type condition for hamiltonicity in tough graphs,
\emph{Electron. J. Combin.}, 29 (2022) Paper 1.5.


\bibitem{Sh} S. Shan, Hamiltonian cycles in $3$-tough $2K_2$-free graphs,
\emph{J. Graph Theory}, 94 (2019) 1--15.

\bibitem{Shan}  S. Shan, Hamiltonian cycles in tough $(P_2\cup P_3)$-free graphs,
\emph{Electron. J. Combin.}, 28 (2021) Paper 1.36.


 \bibitem{SS} L. Shi and S. Shan,
A note on hamiltonian cycles in $4$-tough $(P_2\cup kP_1)$-free graphs,
\emph{Discrete Math.}, 345 (2022) 113081.


\bibitem{Wei} B. Wei,  Hamiltonian cycles in $1$-tough graphs,
\emph{Graphs Combin.}, 12 (1996) 385--395.

\bibitem{West1996} D.B. West, Introduction to Graph Theory, Prentice Hall, Inc., 1996.

\bibitem{Xu} L. Xu, C. Li and B. Zhou, Hamiltonicity of $1$-tough $(P_2\cup kP_1)$-free graphs, \emph{Discrete Math.}, 347(2) (2024) 113755.

\bibitem{Zheng} W. Zheng, H. Broersma and L. Wang,
Toughness, forbidden subgraphs, and Hamilton-connected graphs, \emph{Discuss. Math. Graph Theory}, 42 (2022), no.1, 187--196.









	


\end{thebibliography}
\end{document}